\documentclass[12pt]{article}
\usepackage[utf8]{inputenc}
\usepackage{cmap}
\usepackage{amsmath}
\usepackage{amssymb}
\usepackage{amscd}
\usepackage{xurl}
\usepackage[english]{babel}
\begin{document}
\author{V. K. Beloshapka}

\newcommand{\CC}{\mathbb C}
\newcommand{\RR}{\mathbb R}
\newcommand{\PP}{\mathbb P}
\newcommand{\Ree}{\operatorname{Re}}
\newcommand{\Imm}{\operatorname{Im}}
\newcommand{\rank}{\operatorname{rank}}
\newcommand{\Span}{\operatorname{span}}
\newcommand{\Orb}{\operatorname{Orb}}
\newcommand{\Stab}{\operatorname{Stab}}
\newcommand{\Coker}{\operatorname{Coker}}
\newcommand{\Herm}{\operatorname{Herm}}
\newcommand{\Lie}{\operatorname{Lie}}
\newcommand{\aut}{\mathfrak{aut}}
\newcommand{\Dz}[1]{\frac{\partial}{\partial z_{#1}}}
\newcommand{\Dw}[1]{\frac{\partial}{\partial w_{#1}}}
\newcommand{\Alg}{\mathcal A}
\newcommand{\Cent}{\operatorname{Cent}}
\newcommand{\mm}{\mathfrak m}
\newcommand{\g}{\mathfrak g}

\date{05.08.2026}

\title{\bf Orbits of Exceptional CR Quadrics and Their Variations}

\maketitle

\begin{abstract}
The orbits of exceptional $CR$ quadrics are studied. It is proved that the orbit of the exceptional $(4,4)$ quadric constructed in \cite{VB26-2} is isolated. Using the algebraization procedure described in \cite{VB26-1}, it is shown that the space of orbits of exceptional quadrics of type $(28,28)$ has positive dimension.
\end{abstract}

\footnote{
Faculty of Mechanics and Mathematics, Lomonosov Moscow State University,
Vorobyovy Gory, 119991 Moscow, Russia, vkb@strogino.ru}

{\bf 1. Introduction}

\vspace{3ex}

Nondegenerate $CR$ quadrics and their holomorphic symmetries are among the central objects of $CR$ geometry.
In particular, this is due to the fact that they are the most holomorphically symmetric $CR$ manifolds.

Let us recall the basic definitions \cite{VB02}. Fix two positive integers $1 \leq n$ and $1 \leq k \leq n^2$. Let $H=(H_1,\dots,H_k)$ be a collection of Hermitian matrices of order $n$. We say that this collection is {\it nondegenerate} if the matrices are linearly independent and have no common kernel; that is, if $H_1 \, z=\dots=H_k \, z=0$, then $z=0$. The collection may be interpreted as a vector-valued Hermitian form $<z,\bar{z}>=(H_1 \,z\cdot\bar{z},\dots,H_k \,z\cdot\bar{z})$ (a $k$-Hermitian form). It also determines a $CR$ quadric of type $(n,k)$ ($n$ is the dimension of the complex tangent space and $k$ is the codimension):
$$Q=\{(z,w),\; z\in\mathbf{C}^n,\; w\in\mathbf{C}^k,\; \mathrm{Im}\,w=<z,\bar{z}>\}.$$
Let $\mathrm{Aut}\,Q$ be the local group of holomorphic automorphisms of the germ of $Q$ at the origin, and let $\mathrm{aut}\,Q$ be the corresponding Lie algebra. Its natural grading is defined as follows. Assign weight 1 to the variables $z$, weight 2 to the variables $w$, and the corresponding negative weights to differentiation. Then $\mathrm{aut}\,Q$ becomes a graded Lie algebra of the form
\begin{equation*}
 aut\,Q=g_{-2}+g_{-1}+g_{0}+g_{1}+\dots+g_{d},   \quad  d < \infty,
\end{equation*}
where $g_j$ is the component of weight $j$. The dimensions of $g_{-2}$ and $g_{-1}$ do not depend on $<z,\bar{z}>$ and are equal to $k$ and $2n$, respectively.
By contrast, the dimensions of the nonnegative components of $\mathrm{aut}\,Q$ depend on the Hermitian form $<z,\bar{z}>$.
Examples of quadrics arising in this context exhibited only the values $0$, $1$, and $2$ for $d$. The first example of an exceptional quadric (that is, one with $d>2$), published in \cite{FM}, was of $CR$ type $(4,5)$. In \cite{GM}, an exceptional quadric of type $(6,4)$ was constructed, together with an infinite series of examples of exceptional quadrics of various types.

In the space $\mathcal{H}_{n\,k}$ of $k$-tuples of Hermitian matrices of order $n$, the subset $\mathcal{E}_{n\,k}$ of exceptional collections is a real semialgebraic subset \cite{VB22}. This set is the disjoint union of classes corresponding to holomorphically equivalent quadrics. For quadrics, holomorphic equivalence reduces to a linear equivalence of the form
$z\rightarrow C\,z,\; w\rightarrow\rho\,w$, where $C\in SL(n,\mathbf{C})$ and $\rho\in GL(k,\mathbf{R})$. Clearly, every quadric equivalent to an exceptional one is itself exceptional. Thus it is meaningful to speak of an exceptional orbit.
The set $\mathcal{E}_{n\,k}$ may be empty. The question of whether $\mathcal{E}_{n\,k}$ is empty or nonempty was discussed in detail in \cite{VB26-1}. In particular, an example of an exceptional $(4,4)$ quadric was constructed there. This example is minimal with respect to type, since if $\mathrm{min}(n,k)\leq3$, then the quadric is not exceptional and $\mathcal{E}_{n\,k}=\emptyset$. In the next section, we show that the orbit of this $(4,4)$ quadric is isolated in $\mathcal{E}_{44}$.

\vspace{5ex}
{\bf 2. Isolation of the exceptional $(4,4)$-orbit}

\vspace{3ex}
Let $\mathcal{H}=\mathcal{H}_{44}=\{(H_1, H_2,H_3,H_4)\}$ be the space of quadruples of Hermitian matrices, and let $\mathcal{E}=\mathcal{E}_{44}$ be the subset of $\mathcal{H}$ consisting of the quadruples of Hermitian matrices that define exceptional quadrics. In what follows, $E_{p q}$ denotes the matrix of order four whose entry $M_{p q}=1$, with all other entries equal to zero.
Consider the point $\widetilde H=(H_1,H_2,H_3,H_4) \in \mathcal{H}$,
where
\begin{equation*}
\begin{aligned}
H_1&=E_{12}+E_{21},\\
H_2&=i(E_{21}-E_{12}),\\
H_3&=E_{22},\\
H_4&=E_{14}+E_{41}+E_{23}+E_{32}.
\end{aligned}
\end{equation*}
It defines the quadric $Q\subset\CC^8$
\begin{equation*}
\begin{aligned}
\Imm w_1&=z_1\bar z_2+z_2\bar z_1,\\
\Imm w_2&=iz_1\bar z_2-iz_2\bar z_1,\\
\Imm w_3&=|z_2|^2,\\
\Imm w_4&=z_1\bar z_4+z_4\bar z_1+z_2\bar z_3+z_3\bar z_2.
\end{aligned}
\end{equation*}
The matrices are linearly independent, and $H_4$ is invertible; hence the tuple is nondegenerate.

\vspace{3ex}

The group $G=SL(4,\CC)\times GL(4,\RR)$ acts on $\mathcal H$ by changes of the complex variable $z$ and real changes of basis in the $w$-coordinates. Namely, if $g=(C,\rho) \in G$ and $<z,\bar{z}>$ is the $(4,4)$-Hermitian form defined by the tuple $H=(H_1,H_2,H_3,H_4)$, then $H'=g \cdot H$ defines the form $(\rho)^{-1}<Cz,\overline{Cz}>$. Note that
$$\mathrm{dim}\, \mathcal{H}=4 \times 4^2=64, \quad \mathrm{dim}\, G=(4^2-1) \times 2 + 4^2=46$$
(all dimensions are real).

\vspace{3ex}

{\bf Theorem 1:}
There exists a Euclidean neighborhood $U\subset\mathcal H$ of $\widetilde H$ such that
\begin{equation*}
\mathcal{E} \cap U =\bigl(G\cdot\widetilde H\bigr)\cap U.
\end{equation*}
That is, in a neighborhood of $\widetilde H$, the set $\mathcal{E}$ coincides with the orbit of $\widetilde H$.
Equivalently, the local moduli space of exceptional quadrics at $\widetilde H$ is zero-dimensional.\\
{\it Proof:} The infinitesimal stabilizer for the action of $SL(4,\CC)\times GL(4,\RR)$ is determined by the system
\begin{equation*}
A^*H_j+H_jA
=
\sum_{\ell=1}^4\rho_{j\ell}H_\ell,
\qquad j=1,\ldots,4,   \qquad \mathrm{Tr}(A)=0,
\end{equation*}
where $A\in gl(4,\CC)$ and $\rho\in gl(4,\RR)$. An exact solution of this system gives $\dim  g_0(\widetilde H)=12$. Therefore
$\dim (G\cdot \widetilde H)=46-12=34$.

Choose a real linear complement
\begin{equation*}
\mathcal H
=
T_{\widetilde H}(G\cdot\widetilde H)\oplus N,
\qquad
\dim \, N=64-34=30.
\end{equation*}
Without loss of generality, we may assume that $\varepsilon=(0,0,E_{11},0)\in N$.

An orbit of a real algebraic group is a smooth, locally closed semialgebraic set \cite{BCR}. By choosing a local submanifold of $G$ transverse to the stabilizer and applying the inverse function theorem, we obtain local coordinates ``orbit~$\times$~slice''. Therefore, to prove the theorem, it suffices to show that
\begin{equation*}
\mathcal E\cap(\widetilde H \times N)\cap U
=
\{\widetilde H\}
\end{equation*}
for some sufficiently small neighborhood $U$.

\vspace{3ex}

For a nondegenerate model quadric, a vector field of weight $3$ has the form
\begin{equation*}
X_3=
\left(\;
 d(w,w)+D(w)(z,z), \;\;
 2i\langle z,\overline{d(w,w)}\rangle
\; \right),
\end{equation*}
where $d(w,w)$ is a symmetric bilinear map, while $D(w)(z,z)$ depends linearly on $w$ and is symmetric bilinear in $z$; both $d$ and $D$ take values in $\CC^4$, and
\begin{equation}\label{eq:criterion}
\left\langle D(u)(z,z),\bar z\right\rangle
=
4i\left\langle
 z,\overline{d(\langle z,\bar z\rangle,u)}
\right\rangle, \qquad
\left\langle
 D(\langle z,\bar z\rangle)(z,z),\bar z
\right\rangle=0.
\end{equation}
This system can be written as
$L(H)X=0$. The number of real variables in (\ref{eq:criterion}) is
$(10 \times 4 + 4 \times 10 \times 4 ) \times 2 = 400,$
whereas the number of real linear equations is $(640+800) \times 2 = 2880$. Thus, (\ref{eq:criterion}) is a linear system whose matrix $\mathcal{L}(H) \in M_{2880 \times 400}$ depends linearly on the coordinates of $H$.
Moreover, $K=\mathrm{Ker}(L(H))=g_3(H)$. The matrix $L_0=L(\widetilde H)$ is a highly sparse integer matrix. A computer calculation (Gaussian elimination over $\mathbf{Q}$) gives $\mathrm{rank}(L_0)=394$ and, consequently,
$\mathrm{dim}\, g_3=\mathrm{dim} \, K=400-394=6$.

Set
\begin{align*}
\tau&:=2w_3z_1-(w_1-iw_2)z_2, \\
X_{1}&=
(-iw_1w_3-2w_1z_2^2+w_2w_3+2iw_2z_2^2+4w_3z_1z_2) \Dz{3}
\\[-1mm]
&\quad +2 i w_3^2\Dz{4}+2 w_3\tau\Dw{4},\\
X_{2}&=
(i(iw_1w_3-2w_1z_2^2-w_2w_3+2iw_2z_2^2+4w_3z_1z_2))\Dz{3}
\\[-1mm]
&\quad +2 w_3^2\Dz{4}+2 iw_3\tau\Dw{4},\\
W&:=\Span_{\RR}\{X_{1},X_{2}\}\subset K.
\end{align*}
$2\Ree X_{1}$ and $2\Ree X_{2}$ are elements of $g_3$.

\vspace{3ex}

Consider a family
\[
H(t)=\widetilde H+t\dot H+O(t^2),
\qquad
X(t)=X_0+tX_1+O(t^2),
\]
where $0\ne X_0\in K$, and suppose that $ L(H(t))X(t)=0$.
Extracting the terms of first order in this identity, we obtain
\begin{equation}\label{eq:first-order-equation}
L_0X_1+
D  L|_{\widetilde H}(\dot H)X_0=0.
\end{equation}
Passing to the cokernel of $L_0$ gives rise to a bilinear map
\begin{equation*}
\beta:N\times K\longrightarrow\Coker L_0,
\qquad
\beta(\dot H,X_0)
=
\left[D L|_{\widetilde H}(\dot H)X_0\right].
\end{equation*}
Equation \eqref{eq:first-order-equation} is solvable for $X_1$ if and only if $\beta(\dot H,X_0)=0$.

A direct computation gives
\begin{eqnarray*}
\mbox{if } X_0 \notin W, \mbox{ then} \, \ker\bigl(\beta(\,\cdot\,,X_0):N\to\Coker L_0\bigr)=0;\\
\mbox{if } 0 \ne X_0 \in W, \mbox{ then} \ker\bigl(\beta(\,\cdot\,,X_0)\bigr)=\RR \, \varepsilon.
\end{eqnarray*}

Let $X_0=aX_{1}+bX_{2}$. After deleting the column corresponding to the direction $\varepsilon$, two $29\times29$ minors of the matrix $\beta(\,\cdot\,,X_0)$ take the form
\begin{equation*}
2^{20}3^5a^7(a^2+b^2)^{11},
\qquad
-2^{20}3^5b^7(a^2+b^2)^{11}.
\end{equation*}
Therefore, for $(a,b)\ne(0,0)$, the rank of the restriction to a complement of $\RR\varepsilon$ is $29$. Direct substitution shows that $\varepsilon$ does indeed lie in the kernel. For $X_0$ with a nonzero component along one of the first four basis vector fields, the rank is $30$.

Thus, the only possible nontrivial first-order term transverse to the orbit has the form
\begin{equation}\label{eq:first-order-only}
\dot H=\lambda\varepsilon,
\qquad
X_0=aX_{1}+bX_{2},
\qquad
\lambda\ne0,
\quad (a,b)\ne(0,0).
\end{equation}

Now consider the family
\begin{equation*}
\begin{aligned}
H(t)&=\widetilde H+t\lambda\varepsilon+t^2H_2+O(t^3),\\
X(t)&=X_0+tX_1+t^2X_2+O(t^3),
\end{aligned}
\end{equation*}
where \eqref{eq:first-order-only} holds. The coefficient of $t^2$ in the equation $L(H(t))X(t)=0$ is
\begin{equation}\label{eq:second-order-equation}
\begin{aligned}
0={}&L_0X_2
+D L|_{\widetilde H}(\lambda\varepsilon)X_1
+D  L|_{\widetilde H}(H_2)X_0+\frac12 D^2  L|_{\widetilde H}
(\lambda\varepsilon,\lambda\varepsilon)X_0.
\end{aligned}
\end{equation}
After eliminating all free unknowns from \eqref{eq:second-order-equation}, we obtain
\begin{equation*}
\lambda^2
\begin{pmatrix}
a^2\\ ab\\ b^2
\end{pmatrix}=0.
\end{equation*}
This is a contradiction. Thus, the only transverse direction allowed by the linearization of the exceptionality condition does not extend to a second-order formal deformation. The theorem is proved.

\vspace{5ex}
{\bf 3. Algebraization and the construction of a non-isolated exceptional $(28,28)$-orbit}

\vspace{5ex}

Let $1<t<2$. Define a real commutative associative unital algebra
\[
\Alg_t=\Span_{\RR}\{1,e_1,e_2,e_3,e_4,p,q\}
\]

by the following table of nonzero products in the maximal ideal
$\mm_t=\Span\{e_1,e_2,e_3,e_4,p,q\}$:
\begin{equation}\label{eq:At}
e_1^2=p,\qquad e_2^2=p+q,\qquad
e_3^2=p+tq,\qquad e_4^2=q.
\end{equation}
All products of distinct elements among
$e_1,e_2,e_3,e_4$ are zero, and
\[
p\mm_t=q\mm_t=0.
\]

It follows immediately that the algebra is commutative and associative, and that
\[
\mm_t^2=\Span\{p,q\},\qquad \mm_t^3=0,
\qquad \dim\Alg_t=7.
\]
Multiplication can be written in coordinates. If
$X=\sum_{\alpha=0}^6x_\alpha a_\alpha$,
$Y=\sum_{\alpha=0}^6y_\alpha a_\alpha$, where
$
(a_0,a_1,a_2,a_3,a_4,a_5,a_6)=(1,e_1,e_2,e_3,e_4,p,q),
$
then $XY=\sum r_\alpha a_\alpha$, where
\begin{equation*}
\begin{aligned}
r_0&=x_0y_0,\\
r_\gamma &=x_0y_\gamma+x_\gamma y_0,\qquad 1\leq \gamma \leq4,\\
r_5&=x_0y_5+x_5y_0+x_1y_1+x_2y_2+x_3y_3,\\
r_6&=x_0y_6+x_6y_0+x_2y_2+t x_3y_3+x_4y_4.
\end{aligned}
\end{equation*}
Multiplication defines a symmetric bilinear map
$$\beta_t: \;  V^2=\Span\{e_1,e_2,e_3,e_4\}^2 \rightarrow W=\Span\{p,q\}.$$
Its two coordinate matrices are $M_p=\operatorname{diag}(1,1,1,0)$ and $M_q=\operatorname{diag}(0,1,t,1)$.
The discriminant of this pencil of matrices is the binary quartic
\begin{equation}\label{eq:quartic}
\Delta_t(u,v)=\det(uM_p+vM_q)=uv(u+v)(u+tv).
\end{equation}
An algebra isomorphism induces a linear change of variables in the discriminant. The four roots
of the discriminant on $\mathbb P^1$, namely $(\infty,0,-1,-t)$, have cross-ratio equal to $t$. It depends on the ordering of the roots and is defined up to the following six transformations induced by permutations:
\begin{equation}\label{eq:anharmonic}
t,\quad 1-t,\quad \frac1t,\quad \frac1{1-t},\quad
\frac{t}{t-1},\quad \frac{t-1}{t}.
\end{equation}
Thus $j(t)=256\,\frac{(1-t+t^2)^3}{t^2(1-t)^2}$ is intrinsically defined. For $1<t<2$, only the first element of the list \eqref{eq:anharmonic} lies again in the interval $(1,2)$. We therefore obtain

\vspace{2ex}
{\bf Proposition 2:} If $\Alg_t\simeq\Alg_s,\; 1<s,t<2$, then $t=s$.

\vspace{4ex}

Next, let $\Alg^{\CC}_{t}=\Alg_t\otimes_{\RR}\CC$ be the complexification of $\Alg_{t}$. 
Applying the algebraization procedure of \cite{VB26-1} to the exceptional $(4,4)$-quadric considered above, we obtain a quadric $Q_t$ in $(\Alg^{\CC}_{t})^4 \oplus (\Alg^{\CC}_{t})^4$
\[
Z=(Z_1,\ldots,Z_4)\in (\Alg^{\CC}_{t})^4 \simeq \CC^{28},
\qquad
W=(W_1,\ldots,W_4)\in  (\Alg^{\CC}_{t})^4 \simeq\CC^{28},
\]
\begin{equation}\label{eq:Qt}
\begin{aligned}
\operatorname{Im}W_1&=Z_1\overline Z_2+Z_2\overline Z_1,\\
\operatorname{Im}W_2&=iZ_1\overline Z_2-iZ_2\overline Z_1,\\
\operatorname{Im}W_3&=Z_2\overline Z_2,\\
\operatorname{Im}W_4&=Z_1\overline Z_4+Z_4\overline Z_1
 +Z_2\overline Z_3+Z_3\overline Z_2.
\end{aligned}
\end{equation}
Each equality in \eqref{eq:Qt} is an equality in a seven-dimensional
real algebra and therefore yields seven ordinary real equations.
Thus $Q_t\subset\CC^{56}$ has CR type $(28,28)$.

The polarized form of the original quadric has full image
in $\CC^4$ and trivial common kernel. Thus the algebraized quadric $Q_t$ is nondegenerate.
Algebraization preserves the weight grading and sends every
infinitesimal vector field on the original surface to an $\Alg_t$-family of vector fields
of the same weight.
Therefore,

\vspace{2ex}

{\bf Proposition 3:} The quadric $Q_t$ is exceptional.

\vspace{2ex}

Let $H(Z,\bar{Z})$ be the 4-Hermitian form of the quadric $Q$. Define its {\it centroid} by
$\; \Cent(H)=\{(T,S):\;T\in\operatorname{End}(\CC^4), \; S\in\operatorname{End} (\RR^4)\}$, where
\begin{equation}\label{centr}
 H(T Z,\overline{Z})=S \; H(Z,\overline{Z}) =H(Z,  \overline{T Z})  \text{for all }Z.  
\end{equation}
This is a real algebra, with multiplication given by componentwise matrix multiplication. Under a linear
equivalence of quadrics, their centroids are conjugate; consequently, the isomorphism
type of the algebra $\Cent(Q)$ is an invariant of the orbit of the action of $G=SL(28,\CC)\times GL(28,\RR)$.

Let $E_{ab}^{(m)}$ denote the $m\times m$ matrix whose $(a,b)$-entry is $1$ and whose other entries are zero.
A direct solution of equations (\ref{centr}) for our exceptional $(4,4)$-quadric $Q'$ yields a four-dimensional algebra with basis
$(1,c_1,c_2,c_3)$, where $c_r=(T_r,S_r)$ and
\begin{equation}\label{eq:basecentroid}
\begin{array}{lll}
T_1=E_{31}^{(4)}+E_{42}^{(4)},
&\quad&S_1=E_{41}^{(4)},\\
T_2=iE_{31}^{(4)}-iE_{42}^{(4)},
&&S_2=E_{42}^{(4)},\\
T_3=E_{32}^{(4)},
&&S_3=E_{43}^{(4)}.
\end{array}
\end{equation}
Moreover, $c_\alpha c_\beta=0 \quad(1\leq \alpha,\beta\leq3)$.
Thus,
\begin{equation}\label{eq:C0}
\mathcal C_0:=\Cent(Q')
\simeq\RR\oplus N,
\qquad N=\Span\{c_1,c_2,c_3\},\quad N^2=0.
\end{equation}

\vspace{2ex}
{\bf Proposition 4:} For $t\ne0,1$, the centroid of the algebraized $(28,28)$-quadric $Q^t$ has the form
\begin{equation}\label{eq:centroidtensor}
\Cent(Q^t)\simeq\mathcal C_0\otimes_{\RR}\Alg_t.
\end{equation}
{\it Proof.} The inclusion from right to left is immediate. If $L_a$ denotes multiplication by
$a\in\Alg_t$, then each pair $(T_r,S_r)$ in \eqref{eq:basecentroid}
gives rise to a pair $(T_r\otimes L_a,\ S_r\otimes L_a)\in\Cent(H^t)$.

For the reverse inclusion, write arbitrary $T$ and $S$ as $4\times4$ block matrices with
$7\times7$ blocks relative to the basis $(1,e_1,e_2,e_3,e_4,p,q)$ and substitute the multiplication table
\eqref{eq:At} into (\ref{centr}). Comparing the coefficients of
$1,e_1,e_2,e_3,e_4,p,q$ gives
\begin{equation}\label{eq:centroidgeneral}
(T,S)=\sum_{r=0}^3\sum_{\alpha=0}^6
\lambda_{r\alpha}
(T_r\otimes L_{a_\alpha},S_r\otimes L_{a_\alpha}),
\qquad \lambda_{r\alpha}\in\RR,
\end{equation}
where $(T_0,S_0)=(I_4,I_4)$ is the identity element. During the elimination,
one divides only by $t$ and $t-1$, both of which are nonzero. Consequently,
the solution space has dimension $4\cdot7=28$, and
\eqref{eq:centroidgeneral} is the general solution. Composition of these pairs
corresponds exactly to the multiplication
$(c_r\otimes a)(c_s\otimes b)=c_rc_s\otimes ab$.
This proves \eqref{eq:centroidtensor}.

\vspace{3ex}

Set $D_t:=\mathcal C_0\otimes\Alg_t$ and consider the gradings
\[
\Alg_t=\RR\oplus V\oplus W,
\qquad
\mathcal C_0=\RR\oplus N,
\]

where $\deg V=\deg N=1$, $\deg W=2$. Then
\begin{equation*}
\begin{aligned}
(D_t)_1&=V\oplus N,\\
(D_t)_2&=W\oplus(V\otimes N),\\
(D_t)_3&=W\otimes N.
\end{aligned}
\end{equation*}
This grading coincides with the grading induced by the powers of the
maximal radical of the algebra $\Alg_t$ and is therefore invariant under isomorphisms.

For $x\in(D_t)_1$, let
$m_x:(D_t)_1\to(D_t)_2$ denote multiplication by $x$.
If $0\ne n\in N$, then $\rank m_n=4$. If $0\ne v\in V$, then
\[
\rank m_v=3+\rank\bigl(\beta_t(v,\cdot):V\to W\bigr).
\]

Since the four roots in \eqref{eq:quartic} are distinct,
$\rank\beta_t(v,\cdot)=1$ precisely on the four eigenlines of the pencil,
while away from them this rank is two. Finally, if $v\ne0$ and $n\ne0$, then
$v\otimes n'+v'\otimes n=0$
implies $v'=\lambda v$, $n'=-\lambda n$; hence
$\rank m_{v+n}\geq6$.
Consequently, the invariant projective set
\[
\{[x]\in\mathbb P((D_t)_1):\rank m_x\leq4\}
\]
has a unique two-dimensional linear component, namely $\mathbb P(N)$;
the remaining components are four points. Thus $N$
can be intrinsically recovered from the algebra $D_t$. The ideal
generated by $N$ is $\Alg_t\otimes N$, and
\begin{equation*}
D_t/(\Alg_t\otimes N)\simeq\Alg_t.
\end{equation*}
Therefore, an isomorphism $D_t\simeq D_s$ implies
$\Alg_t\simeq\Alg_s$, and for $1<s,t<2$ it follows that $s=t$.
We have thus proved the following:

\vspace{2ex}

{\bf Theorem 5:} The family $\{Q_t:1<t<2\}$ consists of nondegenerate exceptional
CR quadrics of type $(28,28)$. If $s\ne t$, the quadrics $Q_s$ and $Q_t$ lie
in distinct orbits of the group
$SL(28,\CC)\times GL(28,\RR)$. Consequently, the set of exceptional
$(28,28)$-orbits has at least one continuous modulus.

\vspace{3ex}

If eight-dimensional algebras are used to algebraize the $(4,4)$-quadric, analogous arguments show that
the space of orbits of exceptional $(32,32)$-quadrics has dimension at least six. The paper \cite{P08} studies
parameter spaces of finite-dimensional commutative associative unital algebras and their invariants. Its results imply that seven is the smallest dimension in which continuous families of pairwise nonisomorphic algebras occur (see Proposition~2.1). This motivates the use of seven-dimensional algebras.

\vspace{3ex}

{\bf 4. Example}

\vspace{3ex}

In conclusion, we give an explicit example of two nonisomorphic four-dimensional local algebras such that algebraization of a $(1,1)$-quadric by means of these algebras yields two equivalent $(4,4)$-quadrics. This example explains why, in the construction of Section 3, it is necessary not only to distinguish the algebras $\Alg_t$, but also to prove that $\Alg_t$ can be recovered from the centroid of the algebraized form.

\vspace{3ex}

Consider the two four-dimensional real algebras
\[
\Alg_\pm=\Span_{\RR}\{1,x,y,p\}
\]
with multiplication table $x^2=p,\; y^2=\pm p,\; xy=xp=yp=p^2=0$.
In both algebras, $\mm_\pm=\Span\{x,y,p\},\;
\mm_\pm^2=\Span\{p\},\; \mm_\pm^3=0,
\; \Alg_\pm/\mm_\pm\simeq\RR.
$
Consequently, $\Alg_+$ and $\Alg_-$ are local and hence cannot be decomposed into a direct sum of nonzero unital algebras.
Multiplication induces on $\mm_\pm/\mm_\pm^2$ a quadratic form with values in the one-dimensional space $\mm_\pm^2$:
\[
q_\pm(ax+by)=(a^2\pm b^2)p.
\]
The form $q_+$ is definite, whereas $q_-$ is indefinite. Equivalently,
the projective null cone of the former is empty, while that of the latter consists of two
real points. This property is preserved by algebra isomorphisms, even
if the image of $p$ is multiplied by a nonzero real constant. Thus, $\Alg_+\not\simeq_{\RR}\Alg_-$.

Now consider the hyperquadric and the real algebra
\[
Q_{\rm h}=\{(z,w)\in\CC^2:\operatorname{Im}w=|z|^2\},
\qquad
\mathcal K=\RR[j]/(j^2+1)\simeq\CC.
\]
Thus, $\mathcal K$ is the algebra of complex numbers regarded as a real algebra. If $Z=z_1+jz_2$ and $W=w_1+jw_2$, then the $\mathcal K$-algebraization of $Q_{\rm h}$ is the nondegenerate quadric $\widehat Q$ of type $(2,2)$:
\begin{equation}\label{eq:Qhat}
\begin{aligned}
\operatorname{Im}w_1&=|z_1|^2-|z_2|^2,\\
\operatorname{Im}w_2&=2\operatorname{Re}(z_1\overline z_2).
\end{aligned}
\end{equation}
Successive algebraization is equivalent to algebraization by the tensor
product of the algebras \cite{VB26-1}. Therefore,
\begin{equation}\label{eq:successivealg}
\widehat Q^{\Alg_\pm}
=(Q_{\rm h}^{\mathcal K})^{\Alg_\pm}
\simeq Q_{\rm h}^{\mathcal K\otimes_{\RR}\Alg_\pm}.
\end{equation}

Moreover, there is an explicit isomorphism of real algebras
\[
\Phi:\mathcal K\otimes\Alg_+\longrightarrow
       \mathcal K\otimes\Alg_-
\]
given by
\begin{equation}\label{eq:PhiApm}
\Phi(j)=j,\qquad \Phi(x)=x,\qquad
\Phi(y)=jy,\qquad \Phi(p)=p.
\end{equation}
Indeed, the only nontrivial relation to check is the one involving the change of sign:
\[
\Phi(y)^2=(jy)^2=j^2y^2=(-1)(-p)=p=\Phi(y^2).
\]

In the ordered basis
\[
(1,j,x,jx,y,jy,p,jp)
\]
the isomorphism $\Phi$ is represented by the real matrix
\begin{equation*}
P:(u_1,u_2,u_3,u_4,u_5,u_6,u_7,u_8)
\longmapsto
(u_1,u_2,u_3,u_4,-u_6,u_5,u_7,u_8).
\end{equation*}
Applying $P$ simultaneously to the eight complex $z$-coordinates and the
eight complex $w$-coordinates, we obtain a linear equivalence
\[
\widehat Q^{\Alg_+}\simeq\widehat Q^{\Alg_-}.
\]
Since $\det P=1$, this equivalence belongs to the group
$SL(8,\CC)\times GL(8,\RR)$. Thus,
\begin{equation*}
\Alg_+\not\simeq\Alg_-,\qquad
       \widehat Q^{\Alg_+}\simeq\widehat Q^{\Alg_-}.
\end{equation*}

The reason is that, upon further algebraization, the quadric $\widehat Q=Q_{\rm h}^{\mathcal K}$
detects the algebra
$\mathcal K\otimes\Alg_\pm$, rather than the distinguished real
factor $\Alg_\pm$ within it. The isomorphism \eqref{eq:PhiApm} mixes these two
factors.

\vspace{3ex}

In preparing this article, the author used Maple and ChatGPT.


\begin{thebibliography}{30}

\bibitem{VB26-2} V.K. Beloshapka, On Exceptional CR-Quadrics: Further Developments. \url{https://arxiv.org/pdf/2609.08501}.

\bibitem{VB26-1} V.K. Beloshapka, Associative algebras in CR-geometry, \url{https://arxiv.org/pdf/2603.16779}.

\bibitem{VB02} V.K. Beloshapka, Real submanifolds in complex space: polynomial models, automorphisms, and classification problems // Russian Math. Surveys 57:1 (2002), 1--41.

\bibitem{FM} F. Meylan, A Counterexample to the 2-jet Determination Chern-Moser Theorem in Higher Codimension,
arXiv:2003.11783v1 [math.CV], 26 Mar 2020.

\bibitem{GM} J. Gregorovi$\check{c}$, F. Meylan, Construction of Counterexamples to the
2-jet Determination Chern-Moser Theorem in Higher Codimension, arXiv:2010.10220v1 [math.CV], 20 Oct 2020;
Jan Gregorovi$\check{c}$, Francine Meylan, Construction of Counterexamples to the 2-jet Determination Chern--Moser Theorem in Higher Codimension, Mathematical Research Letters, 29 (2022), no. 2, 399--420. DOI: \url{10.4310/MRL.2022.v29.n2.a4}.

\bibitem{VB22} V.K. Beloshapka, On Exceptional Quadrics, \url{https://arxiv.org/abs/2106.08246}, 15 Jun 2021;
Russian Journal of Mathematical Physics, Vol. 29, No. 1, 2022, pp. 11--27.

\bibitem{BCR}
J.~Bochnak, M.~Coste, M.-F.~Roy,
\textit{Real Algebraic Geometry},
Ergebnisse der Mathematik und ihrer Grenzgebiete, vol.~36,
Springer, Berlin, 1998.

\bibitem{P08}
B.~Poonen,
\textit{Isomorphism types of commutative algebras of finite rank over an
algebraically closed field},
Contemporary Mathematics, vol.~463 (2008), pp.~111--120.


\end{thebibliography}
\end{document}